# ON THE STRENGTH OF SOME TOPOLOGICAL LATTICES

MARCUS TRESSL

– *Dedicated to Murray Marshall* –

ABSTRACT. We study the model theoretic strength of various lattices that occur naturally in topology, like closed (semi-linear or semi-algebraic or convex) sets. The method is based on weak monadic second order logic and sharpens previous results by Grzegorczyk. We also answers a question of Grzegorczyk on the 'algebra of convex sets'.

## Contents

## 1. INTRODUCTION

A lattice in this note is a partially ordered set that has finite suprema and infima. We are concerned with lattices arising naturally in topology, e.g. as zero sets of certain algebraic structures of continuous functions on topological spaces. Examples are:

- Closed subsets of $\mathbb{R}^n$: zero sets of continuous functions $\mathbb{R}^n \longrightarrow \mathbb{R}$.
- Zariski closed subsets of $\mathbb{C}^n$: zero sets of polynomials in $\mathbb{C}[X_1, \ldots, X_n]$.
- Linear subspaces of $K^n$ for a field $K$: zero sets of linear functions $K^n \longrightarrow K$.
- Convex hulls of finite sets in $K^n$ for some ordered field $K$: zero sets of convex continuous semi-linear functions $K^n \longrightarrow K$.

We call lattices of such kind *topological lattices*; there is no formal definition given, but we use the terminology to explain ideas. Given such a lattice $L$, we are interested in the strength of the first order structure $(L, \leq)$, i.e., what is interpretable in the partially ordered set $(L, \leq)$? An important sub-question is: When is $(L, \leq)$ decidable? Many of such decidability questions have been addressed in A. Grzegorczyk's paper on *Undecidability of some topological theories*, [Grzego1951]. We revisit parts of Grzegorczyk's paper, improve some results as well as answer questions from [Grzego1951, end of §5]. I have included footnotes en route when there are strong contact points to [Grzego1951].

On the other hand, our method and our context is different from those in [Grzego1951]. We study topological lattices with the aid of weak monadic second order structures: These are first order structures where it is allowed to quantify over finite subsets of the universe (however we will use a first order implementation of this second order structure, see section 2).

Conceptually, our results say that a topological lattice is decidable, when it comes from a 1-dimensional topological space, and otherwise it interprets $(\mathbb{N}, +, \cdot)$ (or its theory).

Here is an outline of the paper. Our main focus are topological lattices related to some ordered context, although we talk briefly about Zariski closed sets and the $p$-adics in section 7. For this introduction, think of a given topological space $X$, like $\mathbb{R}^n$ or $\mathbb{Q}^n$ and the lattice as a lattice of closed subsets of $X$, like closed semi-linear sets, or closed connected 1-dimensional semi-algebraic subsets of $\mathbb{R}^2$, or the lattice of all closed subsets of $X$.

In many lattices $L$ of closed sets one can find isomorphic (in a sense to be made precise in the text) copies of a unit interval $I$, like the unit interval in $\mathbb{R}$ or in $\mathbb{Q}$. This is explained in section 4. The overall strategy is to analyse which structure is imposed by the ambient $L$ on these intervals. It turns out that $L$ imposes the weak monadic order structure of the total order on $I$. More precisely $L$ defines the partially ordered set of finite unions of closed subintervals $[a, b]$ of $I$. In section 3 we will study this partially ordered set first and show that it is essentially the same as the ordered set $I$ extended by quantification over finite sets (cf. 3.2). Using Rabin's landmark result on the decidability of the monadic second order theory of two successor functions (cf. [Rabin1969]) it follows that this extension is decidable (cf. 3.7).

In particular topological lattices are decidable, when the space itself is 1-dimensional (in an informal sense). If the space is not 1-dimensional then one can frequently find a first order definition of "equal size" of finite sets (these are present as elements in the lattice). This phenomenon suffices to interpret arithmetic in the lattice (cf. 3.3).

In section 6 then, the method leads to interpretations of $(\mathbb{N}, +, \cdot)$ in various topological lattices that are not based on a 1-dimensional space, see 6.2, 6.6, A particular case is given by lattices of convex (not necessarily closed) sets in some $K^n$, $K$ an ordered field. Here we show, using the strategy above, that $L$ indeed interprets the ordered field $K$ expanded by integers. This answers a question of Grzegorczyk at the end of [Grzego1951, §5].

Our standard model theoretic set up follows Hodges book [Hodges1993]. The main reference for o-minimal structures is [vdDries1998]; however the text does not strongly depend on o-minimality and one can just read "semi-linear" or "semi-algebraic" instead of "o-minimal".

## 2. The weak monadic structure of a first order structure

A common feature of the topological lattices that we are concerned with is that they frequently interpret the (weak) monadic second order logic of a structure. We will encompass this logic itself in a first order structure and first look at some model theory of these structures.

**2.1. Definition.** Let $M$ be a first order structure in a language $\mathscr{L}$. The **monadic second order structure of** $M$ is defined to be the following first order structure $\mathrm{MSO}(M)$: The universe is the powerset of $M$. Then $\mathrm{MSO}(M)$ is the expansion of the partially ordered set given by inclusion of subsets of $M$, by the 0-definable subsets of $M^n$; here we identify the set of atoms of $(\mathrm{MSO}(M), \subseteq)$ with $M$ via the bijection $M \longrightarrow \mathrm{Atoms}(\mathrm{MSO}(M))$, $a \mapsto \{a\}$. [1] Notice that the elements of $M$ are, a priori, not 0-definable in $\mathrm{MSO}(M)$.

The **weak monadic second order structure of** $M$ is defined to be the substructure $W(M)$ of $\mathrm{MSO}(M)$ induced on the finite subsets of $M$. Explicitly, $W(M)$ is the first order structure expanding the partially ordered set of finite subsets of $M$ by the 0-definable (in $M$) subsets of $M^n$ (where again we identify $M$ with the atoms of the partially ordered set $W(M)$).

Many undecidability results in this note rely on the undecidability of $W((\mathbb{N},+))$. We address this first, and, for subsequent use, in a slightly more general setting. If $S = (S, \cdot)$ is a semi-group (i.e. $\cdot$ is associative) then we write $X \cdot Y = \{x \cdot y \mid x \in X, y \in Y\}$ for $X, Y \subseteq S$ and $s^K = \{s^k \mid k \in K\}$ for $K \subseteq \mathbb{N}$.[2] An element $s \in S$ is called **torsion** if $s^{\mathbb{N}}$ is finite. Notice that $s^{\mathbb{N}}$ does not need to be a cancellative. On the other hand, if $s$ is not torsion, then $(s^{\mathbb{N}}, \cdot) \cong (\mathbb{N}, +)$.

**2.2. Proposition.** *If $S = (S, \cdot)$ is a semigroup, then the binary relation $t \in s^{\mathbb{N}}$ of $S$ (where $S$ is viewed as a subset of $W(S)$) is 0-definable in $W(S)$, independent of $S$; i.e., the defining formula is the same for all $S$.*

*Proof.* In this proof, "definable" means 0-definable in $W(S)$. We proceed in several steps.

(a) Obviously, the map $W(S) \times W(S) \longrightarrow W(S)$, $(X, Y) \mapsto X \cdot Y$ is definable. Further it is clear that the set of finite semigroups (by which we mean sub-semigroups of $S$) is definable. Consequently the torsion elements of $S$ are definable by saying that $s \in X$ for some finite semigroup $X$.

(b) The following property of $(s, X)$, $s \in S$, $X \in W(S)$ is definable:

$(*)$ $\qquad$ $s$ is not torsion and $\exists n \in \mathbb{N} : X = s^{\{1,\ldots,n\}}$.

*Proof.* It suffices to show that $(*)$ holds if and only if the following properties hold (as "torsion" and the other properties are 0-definable in $W(S)$):

$\qquad s \in X$, $s$ is not torsion and

---

[1] We will in general not set up an explicit language for $MSO(M)$ as we will only talk about definability. Further we choose a single sorted set up, since this fits better with applications later. (Another set up could be two introduce a new sort for the powerset of $M$ and the element relation between $M$ and its powerset.)

[2] $\mathbb{N} = \{1, 2, 3, \ldots\}$.

$$\exists x \in X \left( X \cap x \cdot X = \emptyset \text{ and } s \cdot (X \setminus \{x\}) \subseteq X \right).$$

If $(*)$ holds, then $x = s^n$ has the required properties as $s^{\mathbb{N}}$ is cancellative.

Conversely, suppose $s \in X$ is not torsion and there is some $x \in X$ with $X \cap x \cdot X = \emptyset$ and $s \cdot (X \setminus \{x\}) \subseteq X$. Let $n \in \mathbb{N}$ be maximal with $s^{\{1,\ldots,n\}} \subseteq X$. Since $s^{n+1} \notin X$ and $s \cdot (X \setminus \{x\}) \subseteq X$ we must have $x = s^n$. Consequently, $X \cap s^n \cdot X = \emptyset$. Now take $y \in X$ and let $k \in \{0, \ldots, n-1\}$ be maximal with $s^k \cdot y \in X$ (where $s^0 \cdot y$ stands for $y$). Since $s^n y \notin X$ we have $s^{k+1} y \notin X$ and therefore $s^k y = x$. Hence $s^k y = s^n$, which implies $y = s^{n-k} \in s^{\{1,\ldots,n\}}$ by cancelation. □

(c) We can now define $t \in s^{\mathbb{N}}$ by saying

   $t$ is in the smallest finite semigroup containing $s$, or, there is some $X \in W(S)$ with $t \in X$ such that property $(*)$ holds for $s$ and $X$.

□

**2.3. Corollary.** *Let $S$ be the semigroup $(\mathbb{N}, +)$. Then in $W(S)$ multiplication of natural numbers is 0-definable.*

*Proof.* By 2.2 we can define, without parameters, division of natural numbers by saying that $k \mid n$ if and only if $n$ is in the (additive!) semigroup generated by $k$. It is routine to check that multiplication is 0-definable in $(\mathbb{N}, +, \mid)$. □

**2.4.** For infinite abelian groups $G$ with $p \cdot G = 0$, $p$ prime, the structure $W(G)$ is also undecidable, but for a different reason: Let $T$ be the common theory of all finite dimensional $\mathbb{F}_p$-vector spaces expanded by 5 subspaces. This is known to be undecidable by [Toffal1997, Example 2, p. 246]. Now notice that $G$ is an $\mathbb{F}_p$-vector space and as $G$ is infinite we can talk about all finite dimensional $\mathbb{F}_p$-subspaces (up to isomorphism) in a definable family of $W(G)$. It is then not difficult to effectively construct for each sentence $\varphi$ in the language of $T$ a sentence $\tilde{\varphi}$ in the language of $W(G)$[3], such that

$$T \vdash \varphi \iff W(G) \models \tilde{\varphi}.$$

As $T$ is undecidable, also $W(G)$ must be undecidable. We omit the details here as we do not use it later on. However it is worth mentioning because it has a nice application in our context:

If $M$ is an infinite set (viewed as a first order structure in the empty language), then $W(W(M))$ is undecidable: The reason is that $W(M)$ expands an infinite abelian group of exponent 2, where addition is symmetric difference of finite sets. It should be said that the undecidability of $T$ here is based on the insolvability of the word problem for the class of finite groups. I do not know if $W(W(M))$ interprets Peano arithmetic.

We see from 2.2 and 2.4 that $W(M)$ may become wild for fairly simple structures $M$. There is one notable exception:

**2.5. Proposition.** *Let $T$ be an $\omega$-categorical theory in a countable language $\mathscr{L}$.*
  (i) *If $M \prec N$ are models of $T$, then $W(M) \prec W(N)$.*
 (ii) *For every $M \models T$, $M$ is stably embedded into $W(M)$. In fact, every $S \subseteq M^n$ that is definable in $W(M)$ with parameters $a_1, \ldots, a_k \in W(M)$ is definable in $M$ with parameters from $a_1 \cup \ldots \cup a_k$.*

---

[3]Say the language $\{\subseteq, +\}$, where $+$ is a (partial) binary function symbol

*Proof.* (i) It suffices to do the case when $M$ is countable, since then in general $M$ is the up-directed union of elementary substructures of $N$. We may replace $N$ by a strongly $\aleph_1$-homogeneous elementary extension of $N$. Take a countable elementary substructure $\mathscr{N}$ of $W(N)$. Let $M_0 = N \cap \mathscr{N}$. It is then easy to verify that $M_0 \prec N$ and in fact $W(M_0) = \mathscr{N}$. Since $M, M_0 \prec N$ are countable they are isomorphic and there is a an automorphism of $N$ mapping $M$ to $M_0$. Since $W(M_0) = \mathscr{N} \prec W(N)$ we get $W(M) \prec W(N)$.

(ii) Let $S \subseteq M^n$ be definable in $W(M)$ by a formula $\varphi(\bar{x}, \bar{a})$, $\bar{x}$ and $n$-tuple, $\bar{a} \in W(M)^k$. Let $E$ be the union of the finite sets $a_1, \ldots, a_k$. By (i) we may replace $M$ by a countable elementary substructure containing $E$. Suppose $S$ is not definable in $M$ with parameters from $E$.

<u>Claim.</u> There are $\bar{b}, \bar{c} \in M^n$ such that $\bar{b} \in S$, $\bar{c} \notin S$ and $\mathrm{tp}^M(\bar{b}/E) = \mathrm{tp}^M(\bar{c}/E)$.

*Proof of the claim.* Since $T$ is $\omega$-categorical, $S_n^M(E)$ is finite by Ryll-Nardzewski ([Hodges1993, Theorem 7.3.1]), say $S_n^M(E) = \{q_1, \ldots, q_m\}$ is of size $m$. For each $i \in \{1, \ldots, m\}$ take an $\mathscr{L}(E)$-formulas $\psi_i(\bar{x})$ isolating $q_i$. Since we assume that $S$ is not definable in $M$ with parameters from $E$, $S$ is not a finite union of sets defined by some $\psi_i(\bar{x})$. Hence for some $i$, $S$ and $M^n \setminus S$ hit $\psi_i[M^n]$. Take $\bar{b} \in S \cap \psi_i[M^n]$ and $\bar{c} \in \psi_i[M^n] \setminus S$. Then $\mathrm{tp}^M(\bar{b}/E) = q_i = \mathrm{tp}^M(\bar{c}/E)$, as required. □

Now take $\bar{b}, \bar{c} \in M^n$ as in the claim. Then $\mathrm{tp}^M(\bar{b}/E) = \mathrm{tp}^M(\bar{c}/E)$. Since $M$ is countable and $T$ is $\omega$-categorical, $M$ is saturated and there is an $E$-automorphism $\sigma$ of $M$ with $\sigma(\bar{b}) = \bar{c}$. But then $\sigma$ extends to an $a_1, \ldots, a_k$-automorphism $\bar{\sigma}$ of $W(M)$ and this contradicts the assumption that $S$ is defined in $W(M)$ by $\varphi(\bar{x}, \bar{a})$. □

2.6. *Remark.* Let $M$ be an $\omega$-categorical structure in a relational language $\mathscr{L}$. If the universal theory of $T$ is finitely axiomatisable, then the age of $M$ (i.e. the set of finitely generated substructures of $M$, cf. [Hodges1993, p. 324]) is a definable family in $W(M)$. This applies, for example, to the countable atomless boolean algebra $A$. Notice that $W(A)$ still is undecidable by 2.4.

The following consequence of 2.3 will be used several times later to check if a structure interprets Peano arithmetic, or even $(\mathbb{N}, +, \cdot)$.

2.7. **Proposition.** *Let $M$ be an infinite set and let $E \subseteq W(M) \times W(M)$ be defined by*

$$E(a,b) \iff a \text{ and } b \text{ are disjoint and of the same size.}$$

*Let $(W(M), E)$ be the expansion of $W(M)$ by $E$.*

  (i) *$(W(M), E)$ defines (without parameters) the relation "a and b have the same size" of $W(M)$. If we identify the equivalence classes with $\mathbb{N}_0 = \{0, 1, 2, \ldots\}$, where $n \in \mathbb{N}_0$ stands for the size of a representative, then $(W(M), E)$ interprets addition of natural numbers on $\mathbb{N}_0$ without parameters.* [4]
  (ii) *Suppose $N$ is an expansion of $(W(M), E)$ and $S(x, \bar{y})$ is a formula of $N$ with the following properties:*
     (a) *For each $\bar{b} \in W(M)^{\bar{y}}$ the set $\{\mathrm{size}(a) \mid a \in W(M), N \models S(a, \bar{b})\}$ is finite.*
     (b) *For all $k \in \mathbb{N}$ and all $n_1, \ldots, n_k \in \mathbb{N}$ there is some $\bar{b} \in W(M)^{\bar{y}}$ such that*
         $$\{n_1, \ldots, n_k\} = \{\mathrm{size}(a) \mid a \in W(M), N \models S(a, \bar{b})\}.$$

---

[4] I do not know if $(W(M), E)$ is undecidable. This is also related to question 1 at the end of the paper.

*Then $N$ interprets $(\mathbb{N}, +, \cdot)$ without parameters. The interpretation only depends on $S$ and not on $M$.*

Notation: *If $N$ is an expansion of $(W(M), E)$ for which there is a formula $S(x, \bar{y})$ satisfying (a) and (b) then we say that $N$ **defines finite sequences of integers**.*

*Proof.* (i). Two finite subsets $a, b$ of $M$ have the same size if and only if $E(a \setminus b, b \setminus a)$. This condition is first order in $(W(M), E)$, hence we may define the relation "$a$ and $b$ have the same size" of $W(M)$. We thus may replace $E$ by this relation. We write $[a]$ for the equivalence class of $a$ modulo $E$. Then the relation

$$A([a], [b], [c]) \iff \exists a', b', c' : E(a, a') \;\&\; E(b, b') \;\&\; E(c, c') \;\&\; a' \cap b' = \emptyset \;\&\; a' \cup b' = c'$$

defines in $(W(M), E)$ the graph of addition induced by $(\mathbb{N}_0, +)$ via the bijection

$$W(M)/E \longrightarrow \mathbb{N}_0; \; [a] \mapsto \mathrm{size}(a).$$

(ii). This follows from (i) and 2.3, since (a) and (b) guarantee that $N$ (or better, the expansion of $(W(M), E)$ by the predicate $S$) interprets $W(\mathbb{N}, +)$. $\square$

## 3. Finite unions of closed intervals

**3.1. Outline.** An important lattice that appears in topologies related to orders is the lattice $L$ of finite unions of closed intervals of the form $[a, b]$ for $a, b \in T$, where $T$ is a dense linear order. We will see that the poset $(L, \subseteq)$ is bi-interpretable with a ceratin weak monadic second order structure $W(M)$ of some first order structure $M$ in 3.2. This will be applied later as follows: In topological lattices one can usually interpret $W(M)$ through the lattice $L$; then if $L$ is essentially not based on a 1-dimensional set (the precise formulation can be found in 6.2), one can define the expansion $(W(M), E)$ of $W(M)$ as in 2.7. Condition (ii) of 2.7 for $(W(M), E)$ then comes for free from the order of $T$ (see 3.3) and so by 2.7, $(L, \subseteq)$ will interpret $(\mathbb{N}, +, \cdot)$ (or at least its theory when parameters are involved).

First some notation. If $T = (T, \leq)$ is any totally ordered set, then consider the **betweenness relation** $B$ defined by $T$: $B$ is a ternary relation and $B(x, y, z)$ holds just if $x \leq y \leq z$ or $z \leq y \leq x$. Betweenness relations will show up axiomatically in 4.1 below. Currently there is no need to talk about this axiomatically. We write $[[x, y]] = [x, y] \cup [y, x]$, for the set of all $z \in T$ between $x, y$ and similarly $((x, y))$ for the set of all $z \in T$ properly between $x, y$. Hence $z \in [[x, y]]$ just means $B(x, z, y)$ and this relation is 0-definable in the structure $(T, B)$. However, one needs two parameters to define the order relation in $(T, B)$.

**3.2. Proposition.** *Let $(T, \leq)$ be an infinite totally ordered set with betweenness relation $B$ and let $L = (L, \subseteq)$ be the partially ordered set of finite unions of closed intervals of $T$ of the form $[a, b]$ (which is a distributive lattice with bottom element). Let $B$ be the betweenness relation of $(T, \leq)$ and consider the first order structure $(T, B)$.*

(i) *The poset $L$ is interpretable in $W(T, B)$ without parameters.*
(ii) *If $T$ is densely ordered, then $W(T, B)$ is definable in the poset $L$ without parameters and – using the interpretation from (i) – $L$ and $W(T, B)$ are bi-interpretable.*

*The definition and the interpretation are independent of $T$.*

*Proof.* We write $M = (T, B)$.

(i). For finite sets $E, F, G \in W(M)$ consider the set

$$A_{E,F,G} := \bigcup \{[[e,f]] \mid e \in E, f \in F,$$
$$[[e,f]] \cap E = \{e\}, [[e,f]] \cap F = \{f\}, [[e,f]] \cap G = \emptyset\}.$$

Clearly $A_{E,F,G} \in L$. Conversely, if $U \in L$, then there are elements

$$e_1 \leq f_1 < g_1 < e_2 \leq f_2 < g_2 < \ldots < g_{n-1} < e_n \leq f_n$$

in $T$ with $U = [e_1, f_1] \cup \ldots \cup [e_n, f_n]$. Then $U = A_{E,F,G}$ with $E = \{e_1, \ldots, e_n\}$, $F = \{f_1, \ldots, f_n\}$ and $G = \{g_1, \ldots, g_{n-1}\}$. We can then interpret $L$ in $W(M)$ without parameters as follows: The universe is $W(M)^3$ and the equivalence relation $(E, F, G) \sim (E', F', G') \iff A_{E,F,G} = A_{E',F',G'}$ is clearly 0-definable in $W(M)$; further it is clear that we can interpret $A_{E,F,G} \subseteq A_{E',F',G'}$ in terms of $E, F, G, E', F', G'$ in $W(M)$ without parameters.

(ii). Assume now that $T$ is densely ordered. As a set, $W(M)$ is 0-definable in $L$ as those $A \in L$ for which every definably connected [5] component is an atom of $L$: Since $T$ is densely ordered we can express in $L$ by a $\leq$-formula (remember: $\leq$ in $L$ is $\subseteq$) that a set $C \in L$ is definably connected: we say that it is not the disjoint join of two nonempty sets in $L$; hence we can also talk about definably connected components of $A \in L$.

Having identified the finite subsets of $T$ as a 0-definable subset of $L$ it suffices to define the betweenness relation on the atoms of $L$ that is naturally given by $T$: Given atoms $\{a\}, \{b\}, \{c\}$ of $L$ we say that $\{b\}$ is between $\{a\}$ and $\{c\}$ if $\{b\}$ is included in the smallest definably connected set $A$ with $\{a\}, \{c\} \subseteq A$.

The bi-interpretability is now routine checking and left to the reader. Both methods (i) and (ii) are independent of $T$. □

**3.3. Proposition.** *Let $(T, \leq)$ be an infinite totally ordered set and let $B$ be the betweenness relation defined by $(T, \leq)$. Let $M$ be the structure $(T, B)$.*

*Let $E \subseteq W(M) \times W(M)$ be defined by*

$$E(a, b) \iff a \text{ and } b \text{ are disjoint and of the same size.}$$

(i) *$(W(M), E)$ defines without parameters and independently of $(T, \leq)$, finite sequences of integers in the sense of 2.7.*

(ii) *$(W(M), E)$ interprets $(\mathbb{N}, +, \cdot)$ without parameters and independently of $(T, \leq)$. More precisely, it interprets the natural definition of addition and multiplication on the equivalence classes of $T$ modulo the relation that identifies sets of equal size.*

*Since $(T, B)$ is 0-definable in $(T, \leq)$ items (i) and (ii) are also true for $M = (T, \leq)$.*

*Proof.* Item (ii) follows from 2.7(ii) once we have shown (i). In order to see (i) we may assume by using 2.7(i) that $E$ satisfies

$$E(a, b) \iff a \text{ and } b \text{ are of the same size.}$$

for all $a, b \in W(M)$.

---

[5] definably connected means definably connected in $(T, \leq)$

We call $x, y \in T$ a jump of $b \in W(M)$ if $x, y \in b$, $x \neq y$ and $((x, y)) \cap b = \emptyset$. For $a, b, c \in W(M)$ we define a relation

$$S(a, b, c) \iff \exists x, y \in b \left( ((x, y)) \cap b = \emptyset \ \& \ E(((x, y)) \cap c, a) \right).$$

Hence $S(a, b, c)$ holds just if there is a jump $x, y$ of $b$ such that the size of $((x, y)) \cap c$ is the size of $a$. Since the relation $z \in ((x, y))$ is definable in $M = (T, B)$ by $B(x, z, y) \ \& \ z \neq x, y$, it is clear that $S$ is 0-definable in $(W(M), E)$. Further, for all $b, c \in M$, the set of $E$-equivalence classes of $a \in W(M)$ with $S(a, b, c)$ is finite: It consist of all the $E$-classes of $a \in W(M)$ such that the size of $a$ is the size of $((x, y)) \cap c$ for some jump $x, y$ of $b$. Hence 2.7(ii)(a) holds.

On the other hand, also 2.7(ii)(b) holds: If $n_1, \ldots, n_k \in \mathbb{N}$, then we may obviously choose finite subsets $b, c$ of $T$ such that the cardinalities of the sets $((x, y)) \cap c$ are $n_1, \ldots, n_k$, when $(x, y)$ runs through the jumps of $b$. Hence $S(u, b, c)$ defines a set of representatives of the elements of $\{n_1, \ldots, n_k\}$. □

We turn to decidability of the lattice of finite unions of closed intervals.

**3.4. Theorem.** *Let $S2$ be the binary tree $2^{<\omega}$ together with the two successor functions $\sigma \mapsto \sigma\hat{\ }1$ and $\sigma \mapsto \sigma\hat{\ }0$. Then any expansion of* MSO$(S2)$ *by naming finitely many elements from $W(S2)$ is decidable.*

*Proof.* This is the main result in [Rabin1969], see [Rabin1969, Theorem 1.1] and [Rabin1969, Corollary 1.9] □

3.5. *Remark.* The following binary relations on S2 are 0-definable in $W(S2)$:
  (i) $\sigma \leq \tau$ defined as '$\tau$ extends $\sigma$'. (One can express that $\sigma$ is in the smallest finite subset of $S2$ containing $\tau$, that is closed under immediate predecessors.)
 (ii) $\sigma \preccurlyeq \tau$ defined as '$\sigma$ is to the right of $\tau$' in the natural horizontal order of the binary tree $2^{<\omega}$. Explicitly, this means $\sigma = \tau$, or, $\sigma\hat{\ }1 \leq \tau$, or, $\tau\hat{\ }0 \leq \sigma$, or $\sigma$ and $\tau$ are incomparable for $\leq$ and if $\gamma$ is the infimum of $\sigma$ and $\tau$ for $\leq$, then $\gamma\hat{\ }0 \leq \sigma$.
(iii) It is clear that the poset $(S2, \preccurlyeq)$ is a countable dense linear order without endpoints. Hence $(S2, \preccurlyeq) \cong (\mathbb{Q}, \leq)$ and therefore we may consider $W(S2)$ as an expansion of $W(\mathbb{Q}, \leq)$.

**3.6. Corollary.** *If $T$ is a densely linearly ordered set, then any expansion of $W(T)$ by naming finitely many elements is decidable.* [6]

*Proof.* Since $T$ is $\omega$-categorical we may assume that $T$ is countable, using 2.5(i). Then $T$ is isomorphic to an interval of $(\mathbb{Q}, \leq)$ with endpoints in $\mathbb{Q}$. Since all these intervals are parametrically definable in $(\mathbb{Q}, \leq)$ we may assume that $T = (\mathbb{Q}, \leq)$. By 3.5(iii) it is therefore enough to show that any expansion of $W(S2)$ by naming finitely many elements is decidable. But this follows from 3.4, since $W(S2)$ is 0-definable in MSO$(S2)$: A subset $S$ of $2^{<\omega}$ is infinite if and only if there is a nonempty subset $Y$ of the down set generated by $S$ for $\leq$ such that for all $\sigma \in Y$ either $\sigma\hat{\ }1 \in Y$ or $\sigma\hat{\ }0 \in Y$. □

---

[6]This is generally attributed to Läuchli, see [Laeuch1968], but I was unable to find a precise reference. So below is a proof relying on Rabin's work.

**3.7. Corollary.** *Let $T$ be a dense linear order. Any expansion of the lattice of finite unions of closed intervals of the form $[a,b]$ by finitely many parameters, is decidable.*

*Proof.* By 3.6 and 3.2. □

**3.8.** The lattice $(\mathbb{Q} \times \mathbb{Q}, \leq)$ with component wise order is obviously definable in the totally ordered set $(\mathbb{Q}, \leq)$. However, the weak monadic second order logic of the lattice $(\mathbb{Q} \times \mathbb{Q}, \leq)$ is not interpretable in the weak monadic second order logic of $(\mathbb{Q}, \leq)$; in other words, we can not talk about finite sets of pairs of rational numbers in the first order structure $(W(\mathbb{Q}, \leq))$: We know from 3.6 that $W(\mathbb{Q}, \leq)$ is decidable; however $W(\mathbb{Q} \times \mathbb{Q}, \leq)$ interprets the theory of $(\mathbb{N}, +, \cdot)$, see 6.6 (and 6.5).

## 4. Defining closed and bounded intervals

In order to apply the results from the previous section we need a method to define a totally ordered set $T$ in a topological lattice. The lattice will then frequently define finite unions of closed intervals of $T$.

**4.1. Definition.** Let $X$ be a set and $B$ be a ternary relation on $X$. We say that $B$ is a **bounded betweenness relation** if there are $a, b \in X$ such that

(a) $B(a, x, y)$ and $B(x, y, b)$ define the same total order $\leq$ on $X$ with smallest element $a$ and largest element $b$.
(b) $B$ is the betweenness relation of $\leq$.

We write $[[x, y]] = \{z \in X \mid B(x, z, z)\}$ and $((x, y)) = [[x, y]] \setminus \{x, y\}$. Then $B$ is called **dense** if $((x, y)) \neq \emptyset$ for all $x \neq y$.

Clearly bounded (dense) betweenness relations are first order axiomatisable in the language that has a ternary relation symbol.

**4.2. Notation.** Let $L = (L, \leq, \wedge, \vee, \bot, \top)$ be a lattice with bottom $\bot$ and top $\top$. An atom of $L$ is an element that is minimal among all elements different from $\bot$. We say that an element $x \in L$ is $L$**-connected** (or just connected if $L$ is clear from the context) if there are no $y, z \in L$, with $y \wedge z = \bot \neq y, z$ and $y \vee z = x$.

For example, if $L$ is the poset of linear subspaces of $\mathbb{C}^n$. Then $L$ is a bounded lattice with $a \wedge b = a \cap b$, $a \vee b = a + b$ (sum of linear spaces), $\bot = \{0\}$, $\top = \mathbb{C}^n$. Then the atoms of $L$ are the one dimensional subspaces and these are the only $L$-connected elements of $L$. Notice that $L$ is not distributive.

**4.3. Definition of $I(x)$.** Let $I(x)$ be the formula in the language of posets saying the following:

I1. $x$ is connected, not an atom
I2. The ternary relation of the atoms contained in $x$, defined as "$u$ is in the smallest connected element $\leq x$ containing $v, w$", is a bounded and dense betweenness relation.[7]
I3. For all atoms $u \leq x$ and every $y \leq x$ with $u \not\leq y$ there are atoms $v, w \leq x$ such that $u \in ((v, w))$ and $((v, w))$ contains no atom $\leq y$. Here $((v, w))$ stands for the set of atoms $\leq x$ that are properly between $v, w$ for the betweenness relation from I2.

---

[7]This implicitly implies that for all atoms $v, w \leq x$ there is a smallest connected $b \leq x$ with $v, w \leq b$.

The idea here is that for an element $a$ in a lattice $L$ (and assuming that $L$ is atomic), $I(a)$ says that the set of atoms of $L$ underneath $a$ should be orderable as a dense linear order with endpoints, in such a way that the elements of $L$ are (or, mimic) a basis of closed sets of the order topology. We prefer the formulation with betweenness relations as this avoids the use of constants.

**4.4. Nota Bene.** For $a \in L$ the formula $I(x)$ is true in $L$ at $a$ if and only if it is true in the lattice $\{b \in L \mid b \leq a\}$ at $a$.

**4.5. Proposition.** *Let $X$ be a Hausdorff space that has a countable basis and let $L$ be a lattice of closed subsets of $X$ with the following properties:*

  (a) *If $A \in L$ with $L \models I(A)$, then all subsets of $A$ that are closed in $A$ are in $L$*
  (b) *If $A \subseteq X$ is homeomorphic to $[0,1]_\mathbb{R}$ then all subsets of $A$ that are closed in $A$ are in $L$.*

*(The largest such lattice is the set of all closed subsets of $X$, the smallest such lattice is the set of all closed subsets of $X$ that are contained in a finite union of homeomorphic copies of $[0,1]_\mathbb{R}$.) Then $I(x)$ defines in $L$ the set of all subsets of $X$ that are homeomorphic to the closed unit interval of $\mathbb{R}$.*

*Proof.* By 4.4 and (b), $I(x)$ is satisfied for subsets of $X$ that are homeomorphic to the closed unit interval of $\mathbb{R}$. Conversely take $T \in L$ with $L \models I(T)$. Since $T$ also has a countable basis we may use 4.4 again and assume $X = T$. Now condition (a) says that $L$ is the lattice of all closed subsets of $X$, in particular $L$-connected is just connected in the usual sense of topology. We pick one of the two linear orders $\sqsubseteq$ on $X$ defined by the betweenness relation $B$ asserted to exist by $I(x)$. Let $p, q \in X$ be the smallest and largest element for $\sqsubseteq$. Since $\sqsubseteq$ has $B$ as its betweenness relation we know

$(*)$        for all $a, b \in X$, the set $C(\{a, b\})$ is the interval $[a, b]_\sqsubseteq$,

where $C(\{a, b\})$ stands for the smallest $L$-connected set from $L$ containing $\{a, b\}$. In particular $[a,b]_\sqsubseteq \in L$ is closed in $X$.

<u>Claim.</u> There is a countable infinite subset $Y$ of $X$ that is dense in $X$ and dense in the linear order $(T, \sqsubseteq)$.

*Proof of the claim.* Firstly, $X$ is infinite, since the betweenness relation $B$ is dense. Since $X$ has a countable basis it also has an infinite countable dense subset $Y$ and we only need to see that $Y$ hits every nonempty open interval $(a, b)_\sqsubseteq$ of $(X, \sqsubseteq)$. However, by $(*)$, the set $A := [p, a]_\sqsubseteq \cup [b, q]_\sqsubseteq$ is a closed subset of $X$. Thus, $Y$ hits $(a, b)_\sqsubseteq$.  □

Now pick $Y$ as in the claim and add $p, q$ if necessary. Since the betweenness relation $B$ is dense, $\sqsubseteq$ is a dense total order on $Y$ with endpoints $p, q$. Therefore there is an isomorphism of ordered sets $f : [0, 1]_\mathbb{Q} \longrightarrow (X, \sqsubseteq)$. We extend $f$ to a map $\bar{f} : [0, 1]_\mathbb{R} \longrightarrow X$ as follows: Take $r \in [0, 1]_\mathbb{R} \setminus \mathbb{Q}$ and define

$$X_r = \bigcup_{s \in [0,r] \cap \mathbb{Q}} [p, f(s))_\sqsubseteq \quad \text{and} \quad Y_r = \bigcup_{s \in [r,1] \cap \mathbb{Q}} (f(s), q]_\sqsubseteq.$$

By $(*)$, both $X_r$ and $Y_r$ are open subsets of $X$. Since they are disjoint and nonempty and as $X$ is connected, there is some $z \in X$ that is neither in $X_r$ nor in $Y_r$. On the other hand, the choice of $f$ and the claim imply that there can at most be one such $z$. Hence we may define $\bar{f}(r) = z$.

Using the claim again, it is clear that $\bar{f}$ is an order isomorphism $[0,1]_\mathbb{R} \longrightarrow (X, \sqsubseteq)$ and thus is a homeomorphism of $[0,1]_\mathbb{R}$ and the order topology given by $\sqsubseteq$ on $X$.

It remains to show that the latter topology is the original topology of $X$. Since every closed $\sqsubseteq$-interval is closed in $X$ by $(*)$ we only need to show that every open subset of $X$ is a union of open $\sqsubseteq$-intervals. This is exactly what is expressed in condition I3 of 4.3.

□

4.6. *Remark.* Without condition I3 of 4.3, the formula $I(x)$ would not define homeomorphic copies of $[0,1]_\mathbb{R}$: For example let $X$ be *Smirnov's Deleted Sequence Topology* on $[0,1]$, see [SteSee1995, no. 64, p. 86]: Let $S = \{\frac{1}{n} \mid n \in \mathbb{N}\}$ and let

$$\tau = \{O \setminus A \mid O \subseteq [0,1]_\mathbb{R} \text{ open and } A \subseteq S\}.$$

Then $X$ still has a countable basis, e.g. take all $O \setminus Y$, where $O$ runs through a countable basis of $[0,1]_\mathbb{R}$ and $Y \subseteq S$ is cofinite. Routine checking shows that all intervals of $[0,1]_\mathbb{R}$ are $\tau$-connected and in the lattice of closed subsets of $X$, the element $[0,1]_\mathbb{R}$ satisfies I1 and I2.

**4.7. Proposition.** *Let $R$ be an o-minimal expansion of a dense linear order and let $X \subseteq R^n$ be definable with parameters. Let $L$ be a lattice of parametrically definable subsets of $X$ that are closed in $X$ and suppose $L$ contains all such sets of dimension $\leq 1$. Then $I(x)$ defines in $L$ the set of all subsets of $X$ that are definably homeomorphic to an interval of the form $[a,b]$ of $R$, $a < b$.*

*Proof.* Using 4.4, it is clear that $I(x)$ is satisfied for all $A \in L$ that are definably homeomorphic to some $[a,b] \subseteq R$.

Conversely, by the cell decomposition theorem [vdDries1998, chapter 3, 2.11, p. 52] we see that every $A \in L$ at which $I(x)$ is true must be definably connected and of dimension 1 (if $\dim A \geq 2$, then there are points $x, y \in A$ such that there is no smallest $L$-connected $B \in L$ containing $x, y$). Then $A$ is a finite union of definably homomorphic copies of intervals of $R$ and routine checking shows that $A$ must be definably homomorphic to some $[a,b] \subseteq R$. □

## 5. Linear spaces and convex sets

In [Grzego1951, §5] the "algebra of convexity" is shown to be undecidable. A prototype of such an algebra is the Boolean algebra of subsets of $\mathbb{R}^n$, $n \geq 2$, together with the convex hull operator. Generalizations of this type of closure operators have shown to be undecidable, e.g. see [Davis2013] and [Dornhe1998]. However lattices of convex sets themselves are much less studied from a logic perspective. At the end of [Grzego1951, §5] the question is raised whether the lattice of closed convex subsets of $\mathbb{R}^n$, $n \geq 2$ is undecidable. We confirm this in a strong sense, see 5.2.

**5.1. Reminder on incidence geometry.**
For any field $k$, the lattice $U$ of sub-vector spaces of $k^3$ interprets the field $k$ after naming four parameters. In fact $(U, \leq)$ and $(k, +, \cdot)$ are bi-interpretable with definable parameters (cf. [Hodges1993, Remark 5, section 5.3, p. 215]). This is done via incidence geometry. We recall this briefly. Fix a 2-dimensional subspace $H$ of $k^3$ (hence a maximal element of $U \setminus \{k^3\}$). Let $P$ be the set of all 1-dimensional subspaces of $k^3$ (these are the minimal elements of $U \setminus \{(0)\}$) that are not contained in $H$ and let $L$ be the set of all 2-dimensional subspaces of $k^3$, except $H$. Then the relation $p \subseteq l$ between $p \in P$ and $l \in L$ is isomorphic to the affine incidence geometry defined by $k$ between points of $k^2$ and affine 1-dimensional subspaces of $k^2$. The isomorphism between the incidence geometries is given as follows: We may assume that $H = k \times k \times \{0\}$. Then the map

$$\tau : k^2 \longrightarrow P;\ (a,b) \mapsto (a,b,1) \cdot k \subseteq k^2$$

is a bijection and the induced map between affine 1-dimensional subspaces of $k^2$ and $L$ that maps $A$ to $\tau(A)$ is a bijection as well.

Having the incidence geometries identified we may define the field $k$ (after naming some parameters): As a set it will be a line in the affine geometry, or a 2-dimensional subspace in $k^3$ in the poset of sub-vector spaces of $k^3$:

Addition of $A$ and $C$ on a line with fixed point 0 is defined as follows: Choose a point $B$ not on the line and then use twice the fact that opposite lines in a parallelogram are of the same length. Of course in an arbitrary field this has to be phrased appropriately, but over an ordered field it is exactly this. The pictures show how to construct $A + C$ out of $0, A$ and $C$ step by step:

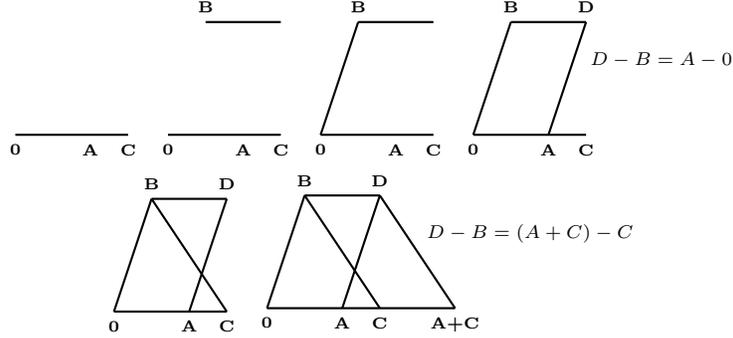

Multiplication of $A$ and $C$ on a line with fixed points 0 and 1 is done by applying the Intercept Theorem twice. Choose a point $B$ not on the line. The pictures show how to construct $A \cdot C$ out of $0, 1, A, C$ step by step:

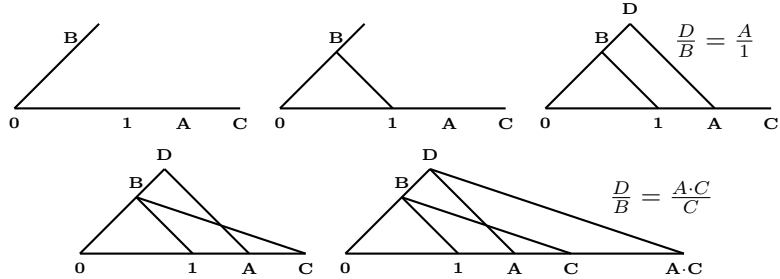

The method of defining "intervals" is particularly easy for convex sets and identifies the strongest topological lattices: these are lattices of convex sets in higher dimensions. For example the lattice of (semi-linear) convex subsets of $\mathbb{R}^n$, $n \geq 2$, interprets $(\mathbb{R}, \mathbb{N}, +, \cdot)$, the expansion of the real field by the set of natural numbers. In detail:

**5.2. Theorem.** *Let $K$ be an ordered field and let $M$ be any expansion of the ordered $K$-vector space $K$. Let $L$ be the poset of all parametrically definable and convex subsets of $M^n$ and let $\overline{L} = \{C \in L \mid C \text{ is closed}\}$. Then the map $L \longrightarrow \overline{L}$* that sends $C$ to its closure is 0-definable in $(L, \subseteq)$; in particular $\overline{L}$ is 0-definable in $L$. The formulas defining these objects are independent of $M$ and $n$.

Further, in the poset $(\overline{L}, \subseteq)$ the following sets are 0-definable:
 (i) The set of all segments $[a, b]$ (= the convex hull of $\{a, b\}$ in $M^n$).
 (ii) The set of affine subspaces of $M^n$.
 (iii) The set of bounded closed convex subsets of $M^n$.
 (iv) If $n \geq 2$, then $\overline{L}$ defines the ordered field $K$ on any 1-dimensional subspace $\ell$ after naming two atoms of $\ell$ (and $\ell$ itself of course).

Now assume $n \geq 2$ and every definable, bounded subset of $M$ that is discrete in $M$, is finite. Then
 (v) The family
$$(\{\ell\} \times F \mid \ell \text{ is a line and } F \text{ is a finite set of atoms} \leq \ell)$$
    is 0-definable in $\overline{L}$.
 (vi) The expansion of the ordered field $K$ by $\mathbb{Q}$ is interpretable in the poset $\overline{L}$.

All formulas in (i)-(vi), defining the various objects can be chosen to be independent of $M$ and $n$.

*Proof.* Firstly, observe that any $A \in L$ is uniquely determined by the atoms $\leq A$. We will talk about points of $A$ rather than atom $\leq A$ when we formulate elementary properties within the poset $(L. \subseteq)$.

<u>Claim 1.</u> The set of all segments $[a, b]$ is 0-definable in $L$.

*Proof.* (Notice that we cannot choose the formula $I(x)$ from 4.3, unless $M = \mathbb{R}$.) A set $C \in L$ is such a segment if and only if $C$ is nonempty, not an atom and the ternary relation of the atoms contained in $C$, defined as "$u$ is in the smallest element of $L$ containing $v, w$", is a bounded and dense betweenness relation. □

We write $I^*(x)$ for the formula in the proof of claim 1.

<u>Claim 2.</u> The map $L \longrightarrow \overline{L}$ that sends $C$ to its closure is 0-definable in $L$.

*Proof.* It is routine to check that indeed for every $C \in L$, the closure of $C$ is in $\overline{L}$ (just copy the proof for $M = \mathbb{R}$). Let $C \in L$. As $C$ is convex, a point $a \in M^n$ is in $\overline{C}$ if and only if for all $b \in M^n$ with $[a, b] \setminus \{a\} \subseteq C$ we have $a \in C$. Hence $\overline{C}$ is the smallest $D \in L$ containing $C$ that has the following property:

> Whenever $S \in L$ with $L \models I^*(S)$ and $A$ is an atom of $L$ with $A \leq S$ such that for all other atoms $B \leq S$ we have $B \leq D$, then also $A \leq D$.

This gives a parameter free definition of the map $L \longrightarrow \overline{L}$ in the poset $L$. □

(i). The formula $I^*(x)$ also defines the segments in $\bar{L} = (\overline{L}, \subseteq)$.

<u>Claim 3.</u> The set of all lines (i.e. 1-dimensional affine subspaces of $M^n$) is 0-definable in $\overline{L}$: A set $A \in \overline{L}$ is a line if and only if it is nonempty, an up-directed union of segments and for all $a \in A$ there is a segment $[b, c] \subseteq A$ with $a \in [b, c] \setminus \{b, c\}$.

(ii). A set $A \in \overline{L}$ is an affine subspace if it is not empty and for all segments $S$ and every line $\ell$ with $S \subseteq A, \ell$ we have $\ell \subseteq A$.

We will from now on assume that $n \geq 2$; item (iii) holds for $n = 1$ by (i).

(iv). By (ii) the affine incidence geometry of $M$ is 0-definable in $\overline{L}$. Hence $(K, +, \cdot)$ is definable in $\overline{L}$ after naming a line $\ell$ and two points $0, 1 \in \ell$ (see 5.1). Furthermore, the order $\leq$ is definable by using the betweenness relation on the segments of $\ell$ and declaring $0 < 1$. This shows (iv).

<u>Claim 4.</u> For each affine hyperplane $A$ of $M^n$ and each line $\ell \nsubseteq A$, let $\pi_{A,\ell}$ be the projection: $M^n \longrightarrow \ell$ along $A$. If we consider $\pi_{A,\ell}$ as a map atoms$(L) \longrightarrow$ {atoms $\leq \ell$}, then the family of (graphs of) all the $\pi_{A,\ell}$ is 0-definable in $\overline{L}$ (by a formula that is independent of $n$ and $M$).

*Proof.* For $a \in M^n$ and $b \in \ell$ we have $\pi(a) = b$ if and only if $a \in \ell$ and $b = a$, or, the unique affine hyperplane containing $a$ and disjoint from $A$, intersects $\ell$ in $b$. This gives a 0-definable definition of the family of all the $\pi_{A,\ell}$. Notice that affine hyper planes are precisely the maximal elements in the poset of proper affine subspaces. □

(iii). A convex set $C$ is bounded if and only if for every projection $\pi_{A,\ell}$ as in claim 4 there is a segment contained in $\ell$ that contains $\pi_{A,\ell}(C)$. (Of course it suffices to test with the $n$ coordinate axes, but in our formulation the definition is independent of $n$.)

For the rest of the proof we now also assume that every definable bounded subset of $M$ that is discrete in $M$, is finite.

<u>Claim 5.</u> The set of all $C \in \overline{L}$ that are convex hulls of finitely many points is 0-definable in $\overline{L}$.

*Proof.* If $F \subseteq K^n$ is finite, then the convex hull $C$ of $F$ in $K^n$ is the convex hull of the extremal points of $C$; this is proved as for $\mathbb{R}^n$. The extremal points are the points $a \in C$ with the property that for all segments $[c, b] \subseteq C$ with $a \in [b, c]$ we have $a = b$ or $a = c$. Further, the extremal points of $C$ are contained in $F$. Hence we can express the property of $C \in \overline{L}$ being such a set as follows:

> $C$ is the smallest element in $\overline{L}$ containing the set $F$ of extremal points (or atoms) of $C$ and any projection of $F$ with a projection function $\pi_{A,\ell}$ as in claim 4 is bounded and discrete in $\ell$.

Using segments we can express "$F$ is discrete in $\ell$" (notice that $F$ is not in $L$); further, by the assumption made right before claim 5, we know that then $F$ must be finite. □

(v). Consider the following property of $A, \ell, C, D \in \overline{L}$:

(∗) $A$ is an affine hyperplane, $\ell \nsubseteq A$ is a line, $C$ is the convex hull of finitely many points and there is an extremal point $E$ of $C$ with $D = \pi_{A,\ell}(E)$

Using claims 4 and 5, this can be expressed by an $\leq$-formula $\varphi(u, v, w, x)$. Consequently, given a line $\ell$ the sets defined by the formulas $\varphi(A, \ell, C, x)$, where $A, C$ vary in $\overline{L}$ is the set of finite subsets (of atoms) of $\ell$.

(vi). By (iv) and (v), the weak monadic second order structure $W(K, +, \cdot)$ is interpretable in $\overline{L}$. By [Bauval1985] this implies that the polynomial ring in one variable over $K$ is interpretable in $\overline{L}$. It is well known that in this ring the rational numbers are definable.

□

6. Lattices of closed definable sets in the ordered context

In this section we fix an ordered field $K$ and work with an expansion $M$ of the ordered $K$-vector space $K$. We will always assume that $M$ is o-minimal, or, $K = \mathbb{R}$ and $M$ is the full structure on $\mathbb{R}$, meaning that all subsets of all $\mathbb{R}^n$ are definable.

6.1. **Lemma.** *Let $n \geq 2$. Let $O \subseteq M^n$ be definable, definably connected, open and of dimension $\geq 2$. If $p_1, q_1, \ldots, p_k, q_k \in O$ are $2k$ points then there are definable homeomorphisms $\sigma_1, \ldots, \sigma_k : [0,1]_K \longrightarrow O$ onto the respective images with mutually disjoint images such that $\sigma_i(0) = p_i$ and $\sigma_i(1) = q_i$ for all $i$.*

*Proof.* This follows by induction on $k$ from the following
<u>Claim.</u> If $\sigma : [0,1]_K \longrightarrow O$ is a definable homeomorphism onto its image, then $O \setminus \sigma([0,1]_K)$ is again definably connected.

Notice that open connected sets are definably path connected under both assumptions on $M$ and that the image of $\sigma$ is closed and definable. The claim then is a routine exercise in patching together paths and left to the reader.  $\square$

6.2. **Theorem.** *Let $n \geq 2$. Let $O \subseteq M^n$ be definable, definably connected, open and of dimension $\geq 2$. Let $L$ be a lattice of definable subsets of $O$ that are closed in $O$ and suppose that $L$ contains all closed subsets of $X$ that are contained in a finite union of definably homeomorphic copies of $[0,1]_K$; when $M$ is the total expansion of $\mathbb{R}$ we also assume that condition (a) of 4.5 is satisfied.[8] Then in $L$ the binary relation*

$$E(A,B) \iff A, B \text{ are finite and of the same size}$$

*is 0-definable in $L$. Further, $(L, \subseteq)$ interprets $(\mathbb{N}, +, \cdot)$ after naming an element $A \in L$ with $L \models I(A)$ (see 4.3).* [9] *The interpretation is independent of $M$, $n$ and $O$.*

*Proof.* Firstly, by 4.5 and 4.7 we can define the property "$A$ is finite" by saying that there is some $C \in L$ with $L \models I(C)$ such that $A$ is a subset of $C$ that is discrete in $C$. By 2.7 it then suffices to show that for disjoint $A, B \in L$ we can define that they are of the same size. By 6.1, this is the case if and only if there is some $C \in L$ such that

(a) each connected component of $C$ hits $A$ in exactly one atom and for each atom $a \leq A$ there is a unique connected component of $C$ containing $a$, and,
(b) each connected component of $C$ hits $B$ in exactly one atom and for each atom $b \leq B$ there is a unique connected component of $C$ containing $b$.

The condition implies that the map from the connected components of $C$ to the atoms of $A$, $D \mapsto D \cap A$ is a bijection, and similarly with $B$. Conversely, if $A$ and $B$ are disjoint and of size $k$, say, then choose $\sigma_1, \ldots, \sigma_k$ according to 6.1. Then the union $C$ of the images of the $\sigma_i$ satisfies (a) and (b).

This shows that $E$ is definable. It follows that for each $A \in L$ with $L \models I(A)$, the structure $(W(A,B), E|_{W(A,B)})$ is definable in $(L, \subseteq)$, where $B$ is the betweenness

---

[8]The largest such lattice is the lattice of all definable subsets that are closed in $O$ and the smallest such lattice is the set of all closed subsets of $O$ that are contained in a finite union of definably homeomorphic copies of $[0,1]_K$.

[9]In [Grzego1951, §3] the undecidability of the lattice of all closed subsets of $\mathbb{R}^n$, $n \geq 2$ is proved; but it is neither shown that the lattice interprets $(\mathbb{N}, +, \cdot)$ nor that the property "finite" is definable in the lattice. See section 8 for a discussion.

relation of $A$ asserted to exist by the formula $I(x)$. Now 3.3 implies that $(\mathbb{N},+,\cdot)$ is interpretable in $L$. □

**6.3. Corollary.** *Let $n \geq 2$ and assume $M$ is o-minimal with the convention of this section in force. Let $P$ be set of all closed definable, definably connected and bounded sets of dimension $\leq 1$, partially ordered by inclusion.*[10] *Then the poset $P$ interprets, independently of $M$, the lattice $L$ of all closed and bounded sets of dimension $\leq 1$. Hence by 6.2, $(P,\subseteq)$ interprets $(\mathbb{N},+,\cdot)$.*

*Proof.* $(L,\subseteq)$ is interpretable in $(P,\subseteq)$ because every set in $L$ is the intersection of two sets in $P$: By o-minimality, there is a definable continuous map $\sigma : [0,1]_K \longrightarrow M^n$ that has $A$ in its image. A suitable perturbation of $\sigma$ will then cut out the complement of $A$ in the image of $\sigma$; the easy details are left to the reader.

Now we can interpret $(L,\subseteq)$ in $(P,\subseteq)$ as $P \times P$ modulo the equivalence relation $(A,B) \sim (C,D) \iff$ the set of atoms below $A$ and $B$ is the set of atoms below $C$ and $D$. Similarly inclusion of $L$ can be interpreted. □

**6.4.** *Remark.* When the dimension of the ambient space is 1 or if the ambient space has no definably connected sets, the lattices in 6.2 generally behave better:
  (i) The lattice of closed subsets of $\mathbb{Q}^n$ is decidable and also the closure algebra of $\mathbb{Q}^n$ is decidable[11]. This is because $\mathbb{Q}^n$ is homeomorphic to $\mathbb{Q}$ and for $\mathbb{Q}$ it is deduced by Rabin from 3.4 in [Rabin1969].
  (ii) In [Rabin1969], Rabin also shows that the lattice of closed subsets of $\mathbb{R}$ is decidable, answering a question of [Grzego1951]. In contrast, Shelah has proved in [Shelah1975] (under the continuum hypothesis) that the closure algebra of $\mathbb{R}$ is undecidable.
  (iii) If $X$ is a Boolean space, then the lattice of closed subsets is decidable, see 7.2(i) below.

**6.5.** It follows from 6.2 that for every ordered field $K$, the lattice of closed and semi-linear subsets of $K^2$ interprets $(\mathbb{N},+,\cdot)$. We can also interpret $(\mathbb{N},+,\cdot)$ in finite unions of rectangles of various sorts. [12] We focus on one case, which also addresses the theme of 3.8: If $T$ is an infinite totally ordered set, then consider the lattice of finite unions of closed infinite rectangles of the form $(-\infty,p] \times (-\infty,q] \subseteq T \times T$ with $p,q \in T$. The partially ordered set of such rectangles is obviously itself a lattice, which is isomorphic to the poset $T \times T$ with componentwise partial order. Further, the elements of $L$ are in bijection with the nonempty, finite anti-chains of $T \times T$ and the order of $L$ translates to the partial order $A \leq B$ of finite anti-chains of $T \times T$ given by $\forall a \in A \exists b \in B : a \leq b$. All these data are 0-definable in $W(T \times T, \leq)$ and so $L$ is 0-definable in $W(T \times T, \leq)$.

**6.6. Proposition.** *Let $T$ be an infinite totally ordered set and let $\leq$ be the componentwise partial order on $T \times T$. Let $L$ be the poset of finite anti-chains of $(T \times T, \leq)$ under the partial order described in 6.5 (hence $L$ indeed is the lattice* described there). Then $L$ defines $(\mathbb{N}, +, \cdot)$ after naming two definable parameters.[13] The definition is independent of $T$.

---

[10] The interest in $P$ comes from the observation that $P$ is the set of images of continuous definable maps $[0,1]_K \longrightarrow M^n$. In a forthcoming paper we will show that the lattice ordered abelian group of these function is decidable, in contrast to the assertion of the proposition.

[11] The closure algebra of $\mathbb{Q}^n$ is the expansion of the boolean algebra of all subsets of $\mathbb{Q}^n$ by the map that sends a set to its closure

[12] In [Grzego1951, §5] the closure operation on the powerset of $\mathbb{R}^2$ with values in such lattices is dealt with, but not the lattices themselves.

*Proof.* In a nutshell, this is true because $L$ is in bijection with the set

$$\{(U, V) \in W(T) \times W(T) \mid U \text{ and } V \text{ have the same size}\};$$

the map sends $A \in L$ to $(p_1(A), p_2(A))$, where $p_i$ are the projections onto the coordinate axes. Hence $L$ "is" the equivalence $E$ of the structure $(W(T), E)$ of 3.3. In detail:

Firstly, the poset $(T \times T, \leq)$ is 0-definable in $(L, \leq)$ as the sub-poset of join irreducible elements. We now define certain data in $(T \times T, \leq)$.

For $x, y \in T^2 = T \times T$ we write $[x, y] = \{z \in T^2 \mid x \leq z \leq y\}$. Let $p, q \in T^2$. We say $p, q$ *define a line segment* if $p \neq q$ and the binary relation $x \leq y$ restricted to $[p, q]$ is a total order with least element $p$ and largest element $q$.

Hence the binary relation "$p, q$ define a line segment" is 0-definable in $(T^2, \leq)$. Obviously, $p, q$ define a line segment if and only if $p \leq q$, $p \neq q$ and ($p_1 = q_1$ or $p_2 = q_2$). A *line segment* is a set of the form $[p, q]$, where $p, q$ define a line segment.

Suppose $p, q \in T^2$ define a line segment. We define

$$\ell(p, q) = \{r \in T^2 \mid [p, r] \cup [r, p] \cup [p, q] \text{ is a line segment}\},$$

which is the line through $p$ and $q$. Hence the ternary relation "$r \in \ell(p, q)$" is 0-definable in $(T^2, \leq)$. Further, we define a map $\pi_{p,q} : T^2 \longrightarrow \ell(p, q)$ by

$$\tau_{p,q}(r) = \begin{cases} r, & \text{if } r \in \ell(p, q), \\ \text{the unique } s \in \ell(p, q) \text{ s.th. } r, s \text{ define a} \\ \text{line segment, or, } s, r \text{ define a line segment,} & \text{if } r \notin \ell(p, q). \end{cases}$$

Hence $\pi_{p,q}$ is the projection onto $\ell(p, q)$.

Let $o, p, q \in T^2$. We say $o, p, q$ *define a coordinate system* if $o, p$ and $o, q$ define line segments and $q \notin \ell(o, p)$, $p \notin \ell(o, q)$. Hence the ternary relation "$o, p, q$ define a coordinate system" is 0-definable in $(T^2, \leq)$.

Now fix $o, p, q \in T^2$ defining a coordinate system. We write $\ell = \ell(o, p)$ and show that the structure $(W(\ell, \leq), E)$ from 3.3 is interpretable in $(L, \leq)$ by using the bijection from the beginning of the proof.

Firstly, we now read the data above in the poset $L$. Hence for example, $\pi_{o,p}$ is a map form the join irreducible elements of $L$ onto the set of join-irreducible elements that correspond to points in $\ell_{o,p}$. Further, we read elements of $L$ as the set of join irreducible elements contained in it. In this set up then, $W(\ell, \leq)$ is interpretable in $L$ by using the identification of $A, B \in L$ when $\pi_{o,p}(A) = \pi_{o,p}(B)$ (which is definable in $L$). Observe that projections are injective on anti-chains.

It remains to interpret "equal size" of elements of $W(\ell, \leq)$: Let $A, B \in L$. Then $\pi_{o,p}(A)$ and $\pi_{o,p}(B)$ have the same size if and only if there are $G, H_A, H_B \in L$ such that
- $\pi_{o,p}(H_A) = \pi_{o,p}(A)$,
- $\pi_{o,p}(H_B) = \pi_{o,p}(B)$,
- $\pi_{o,q}(H_A) = \pi_{o,q}(G)$,
- $\pi_{o,q}(H_B) = \pi_{o,q}(G)$.

---

[13]This means that the set of parameters needed for the definition is 0-definable.

This is first order expressible in $(L, \leq)$, using $o, p, q$ as parameters.

Hence indeed $(W(\ell, \leq), E)$ is interpretable in $L$ and 3.3 entails that $(L, \leq)$ interprets $(\mathbb{N}, +, \cdot)$ after naming the two definable parameters $o, p$; the appearance of $q$ above can be wrapped into an existential quantifier. $\square$

## 7. Zariski closed sets and p-adic sets

**7.1. Zariski closed sets.** Let $n \geq 2$, and let $K$ be an algebraically closed field. Then the lattice $L$ of Zariski closed subsets of $K^n$ interprets $(\mathbb{N}, +, \cdot)$ without parameters and independently of $K$ and $n$.[14]

*Proof.* Firstly, the set of finite subsets of $K^n$ is 0-definable in $L$ as those sets whose irreducible components are atoms. We will check conditions (i),(ii) of 2.7 for $W(K^n)$.

(a) Let $p_1, \ldots, p_k, q_1, \ldots, q_k \in K^n$ be $2k$ points. Then there is a Zariski closed subset $V \subseteq K^n$ such that for each irreducible component $C$ of $V$ there is some (necessarily unique) $i \in \{1, \ldots, k\}$ with $C \cap \{p_1, \ldots, p_k, q_1, \ldots, q_k\} = \{p_i, q_i\}$, and the resulting map from the set of components of $V$ to $\{1, \ldots, k\}$ is a bijection. This property implies that the structure $(W(K^n), E)$ from 2.7(i) is 0-definable in $(L, \subseteq)$.

*Proof.* Firstly, if $E, F \subseteq K^n$ are finite and disjoint. Then there is an irreducible curve $C \subseteq K^n$ containing $E$ and disjoint from $F$: Take an infinite set of irreducible curves passing through $E$ that do not mutually intersect in any other point. Then one of them will do the job.

Using this observation for each $i \in \{1, \ldots, k\}$ and writing $P = \{p_1, \ldots, p_k\}$ and $Q = \{q_1, \ldots, q_k\}$ we get an irreducible Zariski closed set $C_i \subseteq K^n$ with $C_i \cap (P \cup Q) = \{p_i, q_i\}$. Then $V = C_1 \cup \ldots \cup C_k$ has the required property. $\square$

(b) Now we define a family as required in 2.7(ii). Define a ternary relation $S$ on $L$ by

$$S(a, b, c) \iff a, c \text{ are finite and there is an irreducible}$$
$$\text{component } x \text{ of } b \text{ with } a = c \cap b.$$

It is clear that $S$ is 0-definable in $L$. For all $b, c \in L$, the set of cardinalities of $a \in L$ with $S(a, b, c)$ is finite: $c$ has to be finite and the size of $a$ is the size of $c \cap x$ for some connected component $x$ of $b$. Hence 2.7(ii)(a) is satisfied for $S$.

On the other hand, if $S = \{s_1, \ldots, s_k\} \subseteq \mathbb{N}$ is finite of size $k$ then pick a set $b \in L$ with exactly $k$ connected components, each infinite and choose a finite set $a_i$ of size $s_i$ in the $i$-th connected component, disjoint from the other components. Then for $c = a_1 \cup \ldots \cup a_k$, the set of cardinalities of $a \in L$ with $S(a, b, c)$ is $S$. $\square$

**7.2. p-adic sets.** In contrast to the ordered case (cf. 6.2), lattices of closed $p$-adic sets are decidable in all dimensions. This again follows from a result of Rabin in [Rabin1969]: He shows that any Boolean algebra with a second order quantifier over ideals is decidable. In topological terms, this means that the lattice of closed subsets of any Boolean space is decidable. Since the one-point compactification of $\mathbb{Q}_p^n$ is Boolean, this can be used to get decidability of $L$; we leave the details here to the reader.

---

[14]The undecidability of this lattice is also established in [Grzego1951, §3]. See section 8 for a discussion.

## 8. Grzegorczyk's paper

We give a brief outline of [Grzego1951] and talk about the parts that have not been addressed in our note. This hopefully help the reader to study the beautiful results in [Grzego1951] (which are a bit dipped in logic formalism).

The main theme of Grzegorczyk's paper is the creation of a topological version of Robinson arithmetic (see §1]) and then to find models of that. This means, he states an axiom system for partially ordered sets (mainly intended to be applied to topological lattices) such that any lattice that is consistent with these axioms is undecidable. In fact, ordinary Robinson arithmetic is interpreted in the topological set up and then Tarski's theorem on interpretations is invoked. As with ordinary Robinson arithmetic, models of the topological axioms do not need to interpret Peano arithmetic, or even the standard model $(\mathbb{N}, +, \cdot)$.

The models of this topological Robinson arithmetic in [Grzego1951, §2, §4 and §5] are then "closure systems". This means that Boolean algebras are furnished with an operator that behaves in one sense or another like taking closures of arbitrary subsets of a topological space. The set up also encompasses closure operators like $S \longmapsto$ convex hull of $S \subseteq \mathbb{R}^n$, see [Grzego1951, §5]. (For a comparison: In our paper the weaker structures of the images of such operators are studied; recall from 6.4(ii) that the two point of views can differ quite substantially.) The set up in [Grzego1951, §2] is stated as a formal axiom system A1-A6, but this system is infinitary (see A4) and not first order. So the class of closure algebras addressed in §2 is not first order. A more down to earth way of reading the results in §2 is: Every closure algebra satisfying A1-A6 is a model of the topological Robinson arithmetic from §1.

In [Grzego1951, §3], topological lattices are studied in the sense of our paper. More precisely they are studied as a particular class of so-called *Brouwerian algebras*, also know as *(co-)Heyting algebras*. However they are still implicitly defined in terms of the infinitary language from §2 and again one might want to read the theorems in that section as "every lattice satisfying the properties stated at the beginning of §3 satisfies the topological Robinson arithmetic from §1". This should be compared with our results in section 6. For example 6.2 is stronger for the lattices addressed there, because we can also define finite elements in the lattice and interpret $(\mathbb{N}, +, \cdot)$; on the other hand 6.2 is weaker because there are lattices outside the ordered context where [Grzego1951, §3] applies. For example Zariski closed sets, see 7.1.

In [Grzego1951, §4] the topic is "the algebra of bodies". Boolean algebras with a binary predicate "$A, B$ are tangent" are studied. The theme here is closely related to what is called mereotopology (see [AiPrvB2007, Definition 2.5, p. 18]): One should think of the Boolean algebra as the algebra of regular *open* subsets of a topological space and then "being tangent" could be read as "the closures intersect", or as "contact" in mereotopological terminology. Our paper is not addressing this interesting subject, instead we refer to the handbook of spatial logics [AiPrvB2007] (for example see p. 69 there).

We conclude with two questions.
(1) Let $P$ be the set of irreducible Zariski closed subsets of $\mathbb{C}^2$. Does the partially ordered set $(P, \subseteq)$ define the affine subspaces of $\mathbb{C}^2$ (and therefore also $(\mathbb{C}, +, \cdot)$ by 5.1)? Does $(P, \subseteq)$ interpret $(\mathbb{N}, +, \cdot)$?
(2) Let $L$ be the lattice of closed subsets of $\mathbb{R}^2$. Does $L$ interpret the real field $(\mathbb{R}, +, \cdot)$? Note that by 6.2 together with standard coding tricks, we can interpret the field of real algebraic numbers in $L$. There is evidence that in the lattice of closed and semi-algebraic subsets of $\mathbb{R}^2$ the real field is not *definable*; this will be explained in a forthcoming paper.

The University of Manchester, School of Mathematics, Oxford Road, Manchester M13 9PL, UK

Homepage: http://personalpages.manchester.ac.uk/staff/Marcus.Tressl/index.php

*E-mail address*: marcus.tressl@manchester.ac.uk